\documentclass{amsart}

\usepackage{amsaddr}

\usepackage{graphicx} 
\usepackage{tikz, pgfplots}
\title{A new lower bound for deterministic pop-stack-sorting}
\author{Morgan Bauer and Keith Copenhaver}
\address{Mathematics Discipline \\
Natural Sciences Collegium \\
Eckerd College}
\email{copenhkj@eckerd.edu}
\date{July 13, 2023}

\usepackage{amsthm}
\usepackage{amsfonts}
\usepackage{amsmath}

\newtheorem{theorem}{Theorem}
\newtheorem{lemma}[theorem]{Lemma}
\newtheorem{observation}[theorem]{Observation}
\newtheorem{proposition}[theorem]{Proposition}
\newcommand{\Pop}{\mathsf{Pop}}
\newcommand{\ceilfrac}[2]{\left \lceil \frac{#1}{#2} \right \rceil}

\theoremstyle{definition}
\newtheorem{example}[theorem]{Example}

\usepackage[margin=2cm]{geometry}

\begin{document}

\begin{abstract}
    The pop-stack-sorting process is a variation of the stack-sorting process. We consider a deterministic version of this process. We prove a lemma which characterises interior elements of increasing runs after $t$ iterations of the process and provide a new lower bound of $\frac{3}{5}n$ for the number of iterations of the process to fully sort a uniformly randomly chosen permutation of length $n$. 
\end{abstract}

\maketitle

\section{Introduction}
Pop-stack-sorting was introduced by Avis and Newborn in \cite{avis1981pop}. It is a variation of the stack-sorting process first presented by Knuth in {\it The Art of Computer Programming} \cite{knuth1968art}. Pop-stack-sorting has been the subject of much research in recent years (see \cite{asinowski2019pop, asinowski2021flip, choi2022image, claesson2019enumerating, claesson2023counting, defant2021enumeration, defant2022meeting, defant2022pop, defant2022coxeter, hong2022pop}). For our purposes, a {\it permutation} is an ordered list of the elements of $\{1, 2, ..., n\}$, where each element appears once. We write permutations in the one-line notation, that is $12...n$. We consider the deterministic algorithm to which Defant refers as the function $\Pop$ \cite{defant2021fertility}. The function $\Pop$ acts as follows: given a permutation $\sigma$, read the permutation from left to right and if there are any consecutive elements which are in decreasing order, reverse the order of those elements. For example, $\Pop(4\hspace{0.02 in} \underline{71} \hspace{0.02 in} \underline{83} \hspace{0.02 in} 6 \hspace{0.02 in} \underline{952})=4 \hspace{0.02 in}\underline{17}\hspace{0.02 in} \underline{38} \hspace{0.02 in} 6 \hspace{0.02 in} \underline{259}$. Let $\mathcal{D}_n^{\Pop}$ be the average number of iterations of $\Pop$ that it takes to transform a~uniformly randomly chosen permutation of length $n$ to the identity permutation $12...n$. In \cite{defant2021fertility}, Defant conjectured that 
$$\lim_{n \rightarrow \infty} \frac{\mathcal{D}_n^{\Pop}}{n} = 1,$$
and asked, in that paper and at the Banff workshop on Analytic and Probabilistic Combinatorics in November 2022, if there exists $c \in \left( \frac{1}{2}, 1 \right]$ such that

$$\liminf_{n \rightarrow \infty} \frac{\mathcal{D}_n^{\Pop}}n \geq c.$$

It was first proven in \cite{ungar19822n} that $c \leq 1$, and reproven in both \cite{asinowski2021flip, albert2022many} using more direct arguments. The lower bound of $c=\frac{1}{2}$ was not proven by Defant in that paper, but we believe that his conjecture was based on an argument similar to a combination of Proposition \ref{no decreasing runs} and Lemma \ref{lichev lemma}.

In \cite{lichev2022lower}, Lichev gave a lower bound of $c=0.503$. We improve this bound by proving the following theorem.

\begin{theorem} \label{thebigone}
    Let $\mathcal{D}_n^{\Pop}$ be the average number of iterations of $\Pop$ that it takes to transform a uniformly randomly chosen permutation of length $n$ to the identity permutation. Then 
    $$\liminf_{n \rightarrow \infty} \frac{\mathcal{D}_n^{\Pop}}{n} \geq \frac{3}{5}.$$
\end{theorem}


\section{Preliminaries}

For brevity, we adopt the notation $\sigma_t:=\Pop^t(\sigma)$ throughout the paper. Following the standard notation, the position of the element $k$ in the permutation $\sigma_t$ is denoted $\sigma_t^{-1}(k)$. An {\it increasing run} in $\sigma_t$ is a maximal set of elements $\{p_1, p_2, ..., p_k\}$, $k \geq 2$, such that $p_1<p_2<\dots<p_k$ and $\sigma_t^{-1}(p_1)=\sigma_{t}^{-1}(p_2)-1=...=\sigma_{t}^{-1}(p_k)-k+1.$ A {\it decreasing run} is defined similarly. Note that, by this definition, we do not refer to single elements as increasing (or decreasing) runs of length one. An element is in the {\it interior} of an increasing (or decreasing) run if it is an element of an increasing (or decreasing) run, but neither the first nor the last element of that run. We also refer to iterations of $\Pop$ simply as iterations, and the $t$th iteration as iteration $t$.
\\


Defant noted without proof in \cite{defant2021fertility} that in the image of a permutation under $\Pop$, there are no decreasing runs of length more than three. 

\begin{proposition} \label{no decreasing runs}
Let $\sigma$ be a permutation. Then $\sigma_1$ has no decreasing runs of length four or more.
\end{proposition}
\begin{proof}
    Suppose that $\sigma_1$ has four adjacent elements, $a, b, c,$ and $d$, in that order, such that $a>b>c>d.$ Since $\Pop$ reverses decreasing runs, if two elements are in increasing order in $\sigma$, then they cannot be in the same decreasing run and therefore they are also in increasing order in $\sigma_1$. Thus, $b$ lies before $c$ in $\sigma.$ If there are no elements between them in $\sigma$, then $b$ and $c$ are in a decreasing run and are reversed in $\sigma_1$. Hence, there must be some element $e$ in between $b$ and $c$ with either $e>b$ or $e<c$. Assume $e>b$. If there are no elements between $c$ and $d$ in $\sigma$, then they are in the same decreasing run in $\sigma$ and $d$ lies before $c$ in $\sigma_1$. Thus the decreasing run including $e$ and $c$ does not include $d$, so that $e$ is between $c$ and $d$ in $\sigma_1$, a~contradiction. A similar argument can be made to show that if $e<c$ then $e$ must be between $a$ and $b$ in $\sigma_1$.~
\end{proof}


Since all decreasing runs in a pop-stack-sorted permutation are of length two or three, and we are intensely concerned with the manner in which individual elements move, we name each possible way of moving or holding position for an element during an iteration after the first.  \\

We say that $\Pop$ causes two elements to {\it switch} if they are in reverse order but not part of a decreasing run of length three. We say that $\Pop$ causes three elements to {\it pivot} if they are in a decreasing run of length three, where the outermost elements both move and the middle of the three elements holds its position. We call the middle element of a pivot the {\it center} of the pivot. If an element is in the interior of an increasing run, or at the start of an increasing run in position $1$, or the end of an increasing run in position $n$, it has the same position after $\Pop$ is applied, and we say that the element {\it stops}. \\



In any given iteration after the first, there are four ways to move: switch left, switch right, pivot left, pivot right, and two ways to hold position: be the center of a pivot, or stop. \\


Consider the example $\sigma=471836952$ from above. Then we have $\sigma_1=\Pop(4\hspace{0.02 in} \underline{71} \hspace{0.02 in} \underline{83} \hspace{0.02 in} 6 \hspace{0.02 in} \underline{952})=4 \hspace{0.02 in}\underline{17}\hspace{0.02 in} \underline{38} \hspace{0.02 in} 6 \hspace{0.02 in} \underline{259}.$ 
Applying $\Pop$ again, we have 
$\sigma_2=\Pop(\underline{41} \hspace{0.02 in} \underline{73}\hspace{0.02 in} \underline{862}\hspace{0.02 in}5 \hspace{0.02 in} 9)=\underline{14}\hspace{0.02 in} \underline{37} \hspace{0.02 in} \underline {268} \hspace{0.02 in}  5 \hspace{0.02 in} 9.$
During iteration 2, the pair 1 and 4 and the pair 3 and 7 switch and the elements 2, 6, and 8 pivot, with 6 as the center of the pivot. The only elements which hold their position during iteration 2 are 6, 5, and 9; 6 is the center of a pivot and both 5 and 9 stop. Note also that the element 6, which is the center of a~pivot in iteration 2, stops during iteration 1.\\


\section{The pop-stop lemma and proof of Theorem \ref{thebigone}}

We begin with two observations regarding elements which stop. Let $\sigma$ be a permutation of length $n \geq 3$.

\begin{observation} \label{pivot}
If three elements pivot during iteration $t$, $t \geq 2,$ then the center of the pivot stops during iteration $t-1$.
\end{observation}
\begin{proof}

    If $k$ moves or is the center of a pivot during iteration $t-1$, then $k$ is in an increasing run with at least one of its neighbors in $\sigma_{t-1}$, thus $k$ cannot be the center of a pivot during iteration $t$, a contradiction. The only remaining possibility is that $k$ stops during iteration $t-1$.
\end{proof}

\begin{observation} \label{moved left}
    If an element moves to the left during iteration $t$, $t \geq 2$, then during iteration $t-1$ it either moves to the left or stops. The corresponding statement for moving right also holds.
\end{observation}
    
\begin{proof}
If an element moves to the right or is the center of a pivot during iteration $t-1$, then it is larger than its neighbor on its left in $\sigma_{t-1}$ and thus cannot move to the left during iteration $t$. Therefore, the element must move to the left or stop during iteration $t-1$.
\end{proof}


The following lemma, which we refer to as the pop-stop lemma, gives some insight into the circumstances under which an element stops during iteration $t$. Such elements are important because pivots (after the first iteration) are always preceded by a stop. Fix some $1 \leq k \leq n$ in a permutation $\sigma$ of length $n$, where $\sigma^{-1}(k)>k$ so that the sorting process eventually causes $k$ to move to the left. Intuitively, $k$ almost always moves to the left until it becomes a left-to-right maximum unless it encounters an element smaller than $k$ to its left and stops or it is the center of a pivot. However, both cases necessitate stops and therefore there must be smaller elements to the left of $k$, and such elements never move to the right of $k$. If $k$ stops during iteration $t$ where $t$ is large, there must be many elements which are smaller than $k$ to its left, or it must already be a left-to-right maximum. Thus, there exists some $t$ which depends on $k$ but not on $\sigma$ such that $k$ cannot stop during iteration $t$ unless it is a left-to-right maximum.

\begin{lemma}[pop-stop lemma] \label{pop stop}
    Let $\sigma$ be a permutation of length $n$ and let $\sigma_t(i)$ be in the interior of an increasing run. Then $\sigma_t(i)$ is either a left-to-right maximum of $\sigma_t$ or there are at least $\left \lceil \frac{t+1}{2} \right \rceil$ elements smaller than $\sigma_t(i)$ to its left in $\sigma_t$. Similarly, $\sigma_t(i)$ is either a right-to-left-minimum of $\sigma_t$ or there are at least $\left \lceil \frac{t+1}{2} \right \rceil$ elements larger than $\sigma_t(i)$ to its right in $\sigma_t$.

\end{lemma}

\begin{proof}
    We proceed by induction on $t$, the number of iterations. If $t=0$ or $t=1$, $\ceilfrac{t+1}{2}=1$, so the statement is true by virtue of $\sigma_t(i)$ being in the interior of an increasing run. \\



    Suppose an element $\sigma_t(i)$ moves left during iteration $t$. Observation \ref{moved left} implies that $\sigma_{t}(i)$ moves left for $k$ consecutive iterations for some $k \geq 1$ (from iteration $t-k+1$ to iteration $t$, inclusive) and $\sigma_{t}(i)$ either stops during iteration $t-k$ or $k=t$. Let $k^*$ be the number of positions that $\sigma_t(i)$ moved in the previous $k$ iterations, so $k^*=\lvert i-\sigma_{t-k}^{-1}(\sigma_t(i)) \rvert$. It must be the case that $k^* \geq k$. \\

    If $k=t$, then there are at least $k^*>\ceilfrac{k+1}{2}=\ceilfrac{t+1}{2}$ elements larger than $\sigma_{t}(i)$ to its right in $\sigma_{t}$. \\

    If $k<t$, then $\sigma_{t}(i)$ stops during iteration $t-k$ and is in the interior of an increasing run in $\sigma_{t-k-1}$. By the induction hypothesis, $\sigma_{t}(i)$ is a right-to-left minimum in $\sigma_{t-k-1}$ or there are $\ceilfrac{(t-k-1)+1}{2}$ elements larger than $\sigma_t(i)$ to its right in $\sigma_{t-k-1}$. If $\sigma_{t}(i)$ is a right-to-left minimum in $\sigma_{t-k-1}$, it is also a right-to-left minimum in $\sigma_t$. If there are $\ceilfrac{(t-k-1)+1}{2}$ elements larger than $\sigma_t(i)$ to its right in $\sigma_{t-k+1}$, then there are at least 


        $$\left \lceil \frac{(t-k-1)+1}{2} \right \rceil+k^*=\left \lceil \frac{t-k+2k^*}{2} \right \rceil \geq \left \lceil \frac{t+1}{2} \right \rceil$$

        elements larger than $\sigma_{t}(i)$ to its right. \\

    It follows that elements which move left in iteration $t$ are either right-to-left minima or have at least $\left \lceil \frac{t+1}{2} \right \rceil$ larger elements on their right in $\sigma_t$.\\

    A symmetric argument shows that elements which move right are either left-to-right maxima or have at least $\left \lceil \frac{t+1}{2} \right \rceil$ smaller elements on their left in $\sigma_{t}$. \\

    If an element $\sigma_{t}(i)$ stops during iteration $t$, then it is in the interior of an increasing run in $\sigma_{t-1}$. By the induction hypothesis, in both $\sigma_{t-1}$ and $\sigma_t$, $\sigma_{t}(i)$ is either a left-to-right maximum or there are $\left \lceil \frac{t}{2} \right \rceil$ elements smaller than $\sigma_t(i)$ to its left, and it is either a right-to-left minimum or there are $\left \lceil \frac{t}{2} \right \rceil$ elements larger than $\sigma_t(i)$ to its right. \\

    Assume that the element $\sigma_t(i)$ is in the interior of an increasing run in $\sigma_t$. Now we proceed by cases. \\
    
    \begin{enumerate}
        \item[Case 1:] The element $\sigma_{t}(i)$ is the center of a pivot during iteration $t$. \\

        It is not possible for $\sigma_{t}(i)$ to be a left-to-right maximum or right-to-left minimum in $\sigma_{t-1}$ or $\sigma_{t-2}$, since $\sigma_{t-1}(i-1) > \sigma_{t}(i)>\sigma_{t-1}(i+1)$, so we only consider the number of elements. \\
        
        By Observation \ref{pivot}, $\sigma_{t}(i)$ stops during iteration $t-1$ and is in the interior of an increasing run in $\sigma_{t-2}$. By the induction hypothesis, in $\sigma_{t-2}$, there are at least $\ceilfrac{t-1}{2}$ elements smaller than $\sigma_{t}(i)$ to its left. The pivot causes one element smaller than $\sigma_t(i)$ to move from the right of $\sigma_{t}(i)$ to its left, so that, in $\sigma_{t}$, there are at least $\ceilfrac{t-1}{2}+1=\ceilfrac{t+1}{2}$ elements smaller than $\sigma_{t}(i)$ to its left. A symmetric argument shows that there are at least $\ceilfrac{t+1}{2}$ elements larger than $\sigma_{t}(i)$ to its right in $\sigma_t$.\\
        
        \item[Case 2:] The element $\sigma_{t}(i)$ moves right or left during iteration $t$.\\
        
            Suppose $\sigma_{t}(i)$ moves left. Then it must either be a right-to-left minimum or have at least $\left \lceil \frac{t+1}{2} \right \rceil$ elements larger than $\sigma_{t}(i)$ on its right in $\sigma_{t}$. The element $\sigma_{t}(i-1)$ either moves right during iteration $t$, or it is in the interior of an increasing run in $\sigma_{t-1}$. In either case, it is either a left-to-right maximum or there are at least $\left \lceil \frac{t}{2} \right \rceil$ elements smaller than $\sigma_{t}(i-1)$ on its left in $\sigma_{t}$. Thus, $\sigma_{t}(i)$ is either a left-to-right maximum or there are at least $\left \lceil \frac{t}{2} \right \rceil+1 \geq \left \lceil \frac{t+1}{2} \right \rceil$ elements smaller than $\sigma_{t}(i)$ on its left in $\sigma_{t}$. \\

            A symmetric argument applies if $\sigma_t(i)$ moves right. \\
            
        \item[Case 3:] The element $\sigma_{t}(i)$ stops during iteration $t$.\\
        
        Each neighbor of $\sigma_{t}(i)$ either moves towards $\sigma_t(i)$ during iteration $t$ or stops. By the same arguments as in Case 2, in $\sigma_{t}$, $\sigma_{t}(i-1)$ is either a left-to-right maximum or there are at least $\left \lceil \frac{t}{2} \right \rceil$ elements smaller than $\sigma_{t}(i-1)$ on its left. Similarly, $\sigma_{t}(i+1)$ is either a right-to-left minimum or there are at least $\left \lceil \frac{t}{2} \right \rceil$ elements larger than $\sigma_{t}(i+1)$ on its right in $\sigma_{t}$. Therefore, in $\sigma_{t}$, $\sigma_t(i)$ is either a left-to-right maximum or there are at least $\left \lceil \frac{t+1}{2} \right \rceil$ elements smaller than $\sigma_{t}(i)$ on its left and $\sigma_{t}(i)$ is either a right-to-left minimum or there are at least $\left \lceil \frac{t+1}{2} \right \rceil$ elements larger than $\sigma_{t}(i)$ on its right. 


    \end{enumerate}
\end{proof}

In order to be in the center of a pivot (after the initial iteration), an element must first stop. If an element is in position $i$ in $\sigma_{2i-4}$ and it is in the interior of an increasing run, the element has at least $\ceilfrac{2i-4+1}{2}=i-1$ smaller elements to its left, so the pop-stop lemma guarantees that if it stops, then it is a left-to-right maximum. Since the center of a pivot has a larger element to its left, $\sigma_{2i-4}(i)$ is not the center of a pivot, and neither are any other elements which are within $i-1$ positions from either end of $\sigma_{2i-4}$. \\

The most an element can move in any iteration after the first is two positions by pivoting. Thus, for every two iterations, an element can move four positions, however, after those two iterations one more position at the end of the permutation cannot be the center of a pivot. Thus, we split the permutation into fifths; during roughly the first $\frac{2}{5}n$ iterations, an element can move about $\frac{4}{5}n$ positions, but in the following iterations $\frac{1}{5}n$ positions can not be centers of pivots. This can be leveraged to yield the following theorem.

\begin{theorem}\label{bound theorem}
Let $\sigma$ be a permutation of length $n$ and let the element $n-i+1$ (the $i$th largest element) be in the $k$th position after one iteration, that is, $\sigma_1(k)=n-i+1$. Let $t^*$ be the least $t$ such that $\sigma_t$ is the increasing permutation.

If $i \geq 2$, then 
$$t^* \geq 2i-3+\left \lfloor \frac{2(n-k-5i+7)}{5} \right \rfloor+\ceilfrac{n-k-5i+7}{5}.$$

If $i=1$ (so that the element considered is $n$), then

$$t^* \geq \left \lfloor \frac{2(n-k)}{5} \right \rfloor+\ceilfrac{n-k}{5}+1.$$

\end{theorem}
\begin{proof}
    If $i \geq 2$, for the first $2i-4$ iterations, the pop-stop lemma implies that position $n-i+1$ and all positions to its right cannot contain the center of a pivot, which does not contribute to our argument. During those iterations, the element $n-i+1$ moves from position $k$ in $\sigma_1$ to the right at most two positions at a time for the next $2i-3$ iterations, for a total of at most $4i-6$ positions moved after iteration $2i-4$. Thus, $n-i+1$ is in position at most $k+4i-6$ in $\sigma_{2i-4}$. From there, it needs to move through at least another $n-i+1-(k+4i-6)=n-k-5i+7$ elements to reach its final position. The element can pivot for a maximum of $\left \lfloor \frac{2(n-k-5i+7)}{5} \right \rfloor$ iterations and then must switch for at least $\ceilfrac{n-k-5i+7}{5}$ iterations to finish, giving a total of at least $2i-4+\left \lfloor \frac{2(n-k-5i+7)}{5} \right \rfloor+\ceilfrac{n-k-5i+7}{5}$ iterations to sort the permutation. The proof for $i=1$ is similar, but does not require the initial $2i-4$ iterations.
\end{proof} 

We note that this bound could be improved by a single iteration for some values of $n, k,$ and $i$ if one were to consider the possibilities of the distance to travel after $2i-4$ iterations mod 5, but such considerations give at most a change of one in the bound at the cost of a more complicated statement. \\

Lichev proved the following lemma (\cite{lichev2022lower}, Observation 1 and Observation 2).

\begin{lemma}\label{lichev lemma}
In a uniformly randomly chosen permutation $\sigma$ of length $n$, there exists a.a.s.~ an element $n-i+1$ such that after the first iteration $|(n-i+1)-\sigma_1^{-1}(n-i+1)| \geq n-n^{2/3}+1-2\log_2(n).$
\end{lemma}

This lemma combined with Theorem \ref{bound theorem} yields Theorem \ref{thebigone} almost immediately.

\begin{proof}[Proof of Theorem \ref{thebigone}.]
By Lemma \ref{lichev lemma}, there exists a.a.s.~ an element $n-i+1$, $ i<n^{2/3}$, such that after the first iteration, $n-i+1$ is to the left of position $k<n^{2/3}+2\log_2 n$. Combining this with Theorem \ref{bound theorem} yields a lower bound of $\left(\frac{3}{5}-o(1) \right)n$ iterations.

\end{proof}

\section{Future directions}
We note that the bound of $\frac{2}{5}n$ is asymptotically the correct upper bound for the {\it possible} number of pivots to the right which any element can make. We present a family of permutations in which the element $n$ pivots $\frac{2}{5}n-1$ times and reaches the final position after $\frac{3}{5}n$ iterations. 

\begin{example}\label{asymmetric pivot perm}
Let $n=5k$, and let $S_1= \{1, 2, ..., k\}$, $S_2 = \{k+1, ..., 2k\}$, $S_3=\{2k+1, ..., 4k\}$, $S_4=\{4k+1, ..., 5k\}$. Construct a permutation $\sigma$ as follows: take elements from $S_2, S_3,$ and $S_4$ in the pattern $S_4, S_3, S_2, S_3, S_4, S_3, \dots, S_2, S_3$ in decreasing order. Once all of the elements from these sets have been exhausted, append each element of $S_1$ in increasing order at the end. Explicitly, this is the permutation $5k, 4k, 2k, 4k-1, 5k-1, 4k-2, 2k-1, 4k-3, 5k-2, 4k-4, 2k-2, 4k-5, ..., k+1, 2k+1, 1, 2, ..., k.$ Figure 1 shows this construction for $k=4$.

In this case, the permutation takes $\frac{4}{5}n$ iterations to be transformed into the identity; although $n$ is in the last position of $\sigma_{\frac{3}{5}n}$, there are other elements which are not in order in $\sigma_{\frac{3}{5}n}$.
\end{example}

We do not believe this bound is attainable for both the elements $1$ and $n$ in the same permutation. The following construction gives a permutation where each of $1$ and $n$ pivot roughly $\frac{3}{8}n$ times.

\begin{example}\label{pivot perm}
Let $n=4k-1$, and let $S_1= \{1, 2, ..., k\}$, $S_2 = \{k+1, ..., 3k-1\}$, and let $S_3=\{3k, ..., 4k-1\}$. Consider the permutation which consists of elements taken in the pattern $S_3, S_2, S_1, S_2, S_3, S_2, S_1, S_2, \dots, S_3, S_2, S_1,$ in decreasing order. Figure 2 shows this construction for $k=4$. \\

The permutation takes roughly $\frac{3}{4}n$ iterations to be transformed into the identity permutation.
\end{example}

\begin{figure}
\begin{tikzpicture}
\begin{axis}[xtick distance=5,%
scatter/classes={%
    a={draw=black}}]
\addplot[scatter,only marks,%
    scatter src=explicit symbolic]%
table[meta=label] {
x y label
1 20 a
2 16 a
3 8 a
4 15 a
5 19 a
6 14 a
7 7 a
8 13 a
9 18 a
10 12 a
11 6 a
12 11 a
13 17 a
14 10 a
15 5 a
16 9 a
17 1 a
18 2 a
19 3 a
20 4 a
    };
\end{axis}
\end{tikzpicture}
\caption{A permutation using the construction in Example \ref{asymmetric pivot perm} with $k=4$. The element 20 pivots $7=\left(\frac{2}{5}\right)20-1$ times and reaches the final position after $12=\left(\frac{3}{5}\right)20$ iterations. The permutation is sorted after 16 iterations, with the pair 4 and 5 and the pair 10 and 11 being the only elements which are out of order after 15 iterations.}
\end{figure}

\begin{figure}
\begin{tikzpicture}
\begin{axis}[xtick distance=5,%
scatter/classes={%
    a={draw=black}}]
\addplot[scatter,only marks,%
    scatter src=explicit symbolic]%
table[meta=label] {
x y label
1 15 a
2 11 a
3 4 a
4 10 a
5 14 a
6 9 a
7 3 a
8 8 a
9 13 a
10 7 a
11 2 a
12 6 a
13 12 a
14 5 a
15 1 a
    };
\end{axis}
\end{tikzpicture}
\caption{A permutation using the construction in Example \ref{pivot perm} with $k=4$. The elements 1 and 15 both pivot five times and reach their final positions after nine iterations. The permutation is sorted after ten iterations.}
\end{figure}

We believe that the pop-stop lemma can be used in future arguments to improve this lower bound in two ways. \\

First, for every two pivots in an ``early" iteration which can be shown to be unlikely, the pop-stop lemma guarantees an ``extra" pivot which is prevented later. For example, if it were possible to show that the average number of pivots in the first $n/5$ iterations for the element which must travel the furthest was $n/10$, the pop-stop lemma gives a corresponding $n/20$ guaranteed switches at the end. \\

Second, our arguments in the proof of Theorem \ref{bound theorem} use the pop-stop lemma to argue that there are no pivots near the left or right of $\sigma_t$. However, The lemma also implies that elements near the top or the bottom of $\sigma_t$ cannot be the center of pivots. For example, in $\sigma_{2i-4}$, if the element $i$ is in the interior of an increasing run, it has at least $\ceilfrac{2i-3}{2}=i-1$ smaller elements to its left, so it cannot be the center of a pivot in $\sigma_{2i-3}$ and beyond, as a pivot requires that the element $i$ must stop in iteration $t$ and have an element smaller immediately to its right in $\sigma_t$. One can picture that the pop-stop lemma guarantees a square where no centers of pivots can occur that starts one element from the edge of the permutation (for the elements in the first and last positions and the elements $1$ and $n$ can never be the center of a pivot) and closes in one element in each direction every two iterations (see Figure 3). As Example 8 shows, it is not {\it impossible} for elements to pivot every time in the first $\frac{2}{5}n-1$ iterations, but using these closing windows at the top and bottom might be one approach to show that it is rare in a uniformly randomly chosen permutation.


\section{Acknowledgements}
The authors are grateful to Colin Defant for introducing us to the problem and for some helpful suggestions and corrections. We are also grateful to Kelly Debure for lending us some time.

\begin{figure}
\begin{tikzpicture}

\begin{axis}[xtick distance=20,%
scatter/classes={%
    a={draw=black}}]
    \draw[color=black, dashed] 
        (axis cs:7.5, 7.5) -- (axis cs:7.5, 92.5);
    \draw[color=black, dashed] 
        (axis cs:92.5, 92.5) -- (axis cs:7.5, 92.5);
    \draw[color=black, dashed] 
        (axis cs:92.5, 92.5) -- (axis cs:92.5, 7.5);
    \draw[color=black, dashed] 
        (axis cs:7.5, 7.5) -- (axis cs:92.5, 7.5);
\addplot[scatter,only marks,%
    scatter src=explicit symbolic]%
table[meta=label] {
x y label
1   15   a
2   19   a
3   25   a
4   33   a
5   41   a
6   17   a
7   47   a
8   52   a
9   58   a
10   67   a
11   68   a
12   71   a
13   48   a
14   85   a
15   22   a
16   73   a
17   32   a
18   93   a
19   3   a
20   99   a
21   6   a
22   76   a
23   8   a
24   78   a
25   21   a
26   79   a
27   1   a
28   82   a
29   18   a
30   98   a
31   51   a
32   75   a
33   38   a
34   55   a
35   9   a
36   57   a
37   16   a
38   62   a
39   13   a
40   65   a
41   20   a
42   63   a
43   80   a
44   45   a
45   26   a
46   56   a
47   95   a
48   28   a
49   72   a
50   36   a
51   5   a
52   90   a
53   10   a
54   84   a
55   37   a
56   42   a
57   12   a
58   100   a
59   27   a
60   40   a
61   91   a
62   39   a
63   44   a
64   46   a
65   49   a
66   11   a
67   69   a
68   2   a
69   60   a
70   92   a
71   30   a
72   94   a
73   14   a
74   50   a
75   61   a
76   59   a
77   23   a
78   64   a
79   31   a
80   86   a
81   7   a
82   89   a
83   24   a
84   70   a
85   35   a
86   87   a
87   4   a
88   74   a
89   29   a
90   43   a
91   53   a
92   34   a
93   54   a
94   66   a
95   77   a
96   81   a
97   83   a
98   88   a
99   96   a
100   97   a
    };
\end{axis}
\end{tikzpicture}
\begin{tikzpicture}
\begin{axis}[xtick distance=20,%
scatter/classes={%
    a={draw=black}}]
        \draw[color=black, dashed] 
        (axis cs:12.5, 12.5) -- (axis cs:12.5, 87.5);
    \draw[color=black, dashed] 
        (axis cs:87.5, 87.5) -- (axis cs:12.5, 87.5);
    \draw[color=black, dashed] 
        (axis cs:87.5, 87.5) -- (axis cs:87.5, 12.5);
    \draw[color=black, dashed] 
        (axis cs:12.5, 12.5) -- (axis cs:87.5, 12.5);
\addplot[scatter,only marks,%
    scatter src=explicit symbolic]%
table[meta=label] {
x y label
1   15   a
2   17   a
3   19   a
4   25   a
5   22   a
6   33   a
7   32   a
8   41   a
9   3   a
10   47   a
11   6   a
12   48   a
13   8   a
14   52   a
15   21   a
16   58   a
17   1   a
18   67   a
19   18   a
20   68   a
21   51   a
22   71   a
23   38   a
24   85   a
25   9   a
26   73   a
27   16   a
28   93   a
29   13   a
30   99   a
31   20   a
32   76   a
33   26   a
34   78   a
35   28   a
36   79   a
37   5   a
38   82   a
39   55   a
40   10   a
41   98   a
42   36   a
43   75   a
44   45   a
45   12   a
46   56   a
47   57   a
48   27   a
49   62   a
50   37   a
51   63   a
52   39   a
53   40   a
54   65   a
55   11   a
56   80   a
57   2   a
58   95   a
59   30   a
60   42   a
61   72   a
62   14   a
63   90   a
64   23   a
65   84   a
66   44   a
67   31   a
68   46   a
69   7   a
70   100   a
71   24   a
72   91   a
73   49   a
74   35   a
75   50   a
76   59   a
77   4   a
78   60   a
79   29   a
80   69   a
81   34   a
82   92   a
83   43   a
84   94   a
85   53   a
86   61   a
87   54   a
88   64   a
89   66   a
90   86   a
91   70   a
92   89   a
93   74   a
94   77   a
95   81   a
96   87   a
97   83   a
98   88   a
99   96   a
100   97   a
    };
\end{axis}
\end{tikzpicture}
\caption{A permutation of length 100 after 10 sorts, with a box excluding the outermost $\left \lceil \frac{10+1}{2} \right \rceil +1=7$ elements in each cardinal direction, and then after 20 sorts, with a box excluding the outermost $\left \lceil \frac{20+1}{2} \right \rceil +1=12$ elements. The elements outside of the box are either left-to-right maxima, right-to-left minima, or they cannot stop until they become one of the two. Therefore no elements outside of the box will be the center of a pivot in any later iteration.}
\end{figure}

\bibliographystyle{abbrv} 
\bibliography{refs} 
\nocite{*}
\end{document}